\documentclass[12pt]{article}
\usepackage{amsmath}
\usepackage{amssymb}
\usepackage{amsthm}

\numberwithin{equation}{section}
\newtheorem{theorem}{Theorem}[section]

\newtheorem{corollary}[theorem]{Corollary}

\newtheorem{definition}[theorem]{Definition}

\newtheorem{lemma}[theorem]{Lemma}

\newtheorem{proposition}[theorem]{Proposition}

\title{Polyhedrality, Complementarity and Regularity\\ with Application to Variational Inequalities over Polyhedral Sets}

\author{Alexander D. Ioffe\thanks{Department of Mathematics, Technion, Haifa 32000;\
        email: \ alexander.ioffe38@gmail.com}}

\def\xyb{(\bar x,\bar y)}

\def\dis{\displaystyle}

\def\ran{\rangle}
\def\lan{\langle}

\def\Pi{{\cal P}}

\def\N{\mathbb{N}}
\def\R{I\!\!R}
\def\la{\lambda}
\def\al{\alpha}
\def\del{\delta}
\def\ep{\varepsilon}

\def\xb{\overline x}

\def\yb{\overline y}

\def\vb{\overline v}
\def\zb{\overline z}

\def\gr{{\rm Graph}~}

\def\rra{\rightrightarrows}

\def\ep{\varepsilon}
\def\la{\lambda}

\def\al{\alpha}

\def\del{\delta}

\def\cf{{\mathcal F}}

\def\cone{{\rm cone}}

\def\intr{{\rm int}~}
\def\im{{\rm Im}~}

\def\conv{{\rm conv}~}

\def\ri{{\rm ri}~}
\def\ran{\rangle}
\def\lan{\langle}

\def\sur{{\rm sur}}

\def\reg{{\rm reg}}

\def\lip{{\rm lip}}

\def\bd{{\rm bd}}

\def\dis{\displaystyle}

\begin{document}

\maketitle

\section{Introduction}
The starting point for this study were the results on uniqueness and Lipschitz behavior
of solutions to linear variational inequalities
\begin{equation}\label{lvi}
z\in Ax+N(C,x),
\end{equation}
with $A$ being a linear operator in $\R^n$ and $C\subset \R^n$ a polyhedral set,
obtained in the 90s by Robinson \cite{SMR92}, Ralph \cite{DaRa93} and Dontchev and Rockafellar \cite{DR96}. Here are the statements (terminology will be explained later):

\begin{theorem}[Robinson (1992)]\label{rob}
The solution map $S(z)$ of (\ref{lvi}) is single-valued and Lipschitz on $\R^n$
if and only if the coherent orientation condition is satisfied.
\end{theorem}
\begin{theorem}[Dontchev-Rockafellar (1996)]\label{dr96}
Assume that $\xb$ is a solution of (\ref{lvi}) corresponding to $\zb$. Then

(a) the set-valued mapping $\Phi(x) = Ax+N(C,x)$ is strongly regular near $(\xb,\zb)$,
provided it is metrically regular near the point;

(b) the solution mapping of (\ref{lvi}) is single-valued and Lipschitz around $(\xb,\zb)$ if and only if the critical face condition is satisfied.

\end{theorem}

These results have long become classics of variational analysis, in particular because of their key role for analysis of stability of solutions in mathematical programming.  It should be mentioned at this point that  the object of primary interest in the quoted papers was a nonlinear version of (\ref{lvi}):
\begin{equation}\label{nlvi}
z\in f(x)+N(C,x)
\end{equation}
with a continuously differentiable $f: \R^n\to\R^n$. However,
in today's variational analysis, equivalence of the problems,
as far as
local regularity properties of (\ref{lvi}) and (\ref{nlvi})
are concerned, is rather a straightforward consequence of one of the main
principles of the regularity theory of variational analysis, the perturbation theorem,
going back to Milyutin (see \cite{DMO} and  \cite{AI15} for further bibliographic details).
Originally,
a strong regularity version of the perturbation theorem specifically designed for (\ref{nlvi})-like mappings was introduced
by Robinson in his 1980 pioneering study of generalized equations \cite{SMR80}
 to justify the reduction to a linear variational inequality.

Surprisingly,
the available proofs of the results (with the exception of part (b) of the second theorem),  use very little  of variational analysis and
are based on a variety of non-elementary facts of topology and algebra of matrices.
To emphasize the point,
Dontchev and Rockafellar write in the second edition of their monograph
\cite{DR}
after stating a theorem (Theorem 4H.1) containing both statements of Theorem \ref{dr96}:  "The most important part of this theorem that, for the
mapping $f +N(C,\cdot)$ with a smooth $f$ and a polyhedral convex set $C$, metric regularity
is equivalent to strong metric regularity, will not be proved here in full generality.
{\it To prove this fact we need tools that go beyond the scope of this book}" (italicized by me, not by the authors).

This is a definite indication of a highly non-trivial nature of the problem. At
the same time it makes it very reasonable, even necessary  to ask  whether it is possible to develop an adequate alternative theory based primarily
on techniques more traditional for variational analysis. This is not just a question of methodology or aesthetics, it is also and mainly about consistency of variational analysis as such: whether it has enough strength to solve problems simply and naturally formulated in its terms.

Another reason to ask the question is that for some special, and probably the
most important cases (the standard complementarity problem or  variational inequalities associated with the KKT conditions)  proofs of equivalence of metric regularity and strong regularity  by means of variational analysis are well  known since
the 90s \cite{DR96} and \cite{BK97} and can be found in monographs of Dontchev-Rockafellar  \cite{DR} and Klatte and Kummer \cite{KK02})\footnote{Note also that there are variational inequalities with non-polyhedral $C$ for which regularity and strong regularity are equivalent, namely when $C$ is an arbitrary convex set and $A$
is symmetric positive semi-definite - see \cite{DH94}.}.

 The purpose of this paper is to give a positive answer  to the question.
 Our approach differs from the approach in \cite{SMR92,DaRa93} which is based
 on reduction of the problem to the "normal map"
 $$
 A_C= (A\circ\Pi_C)(y)+ (I-\Pi_C),
 $$
 where $\Pi_C$ is the projection  onto $C$ (that is $\Pi_C(y)$ is the closest to $y$ point of $C$), and subsequent use of well developed theory of piecewise affine mappings along with
 a number of non-trivial algebraic and topological arguments.
  Instead, we consider a broader class of set-valued mappings involving certain ``complementarity relation" on a polyhedral set. If $C$ is a cone,
 the relation reduces to a complementarity correspondence between faces of two polyhedral cones extending the well known connection  between faces of a polyhedral cone and its polar. It is not clear whether there are other interesting applications that can be covered by this approach but  simplification for understanding of variational inequalities over polyhedral sets it brings about seems to be really dramatic.  It would be wrong to claim that our proofs are elementary. They are rather not. But they do not use anything beyond basic facts from geometry of polyhedral sets and metric regularity
 theory. Otherwise the paper is fully self-contained.

The plan of the paper is the following. The two subsequent short sections contain  necessary information about metric regularity and polyhedral sets. In \S~4 we introduce the mentioned concepts of complementarity relation on a polyhedral set and   of face complementarity of polyhedral cones. Section 5 is central.
We prove there counterparts of Theorems \ref{rob} and \ref{dr96} for a class of set-valued mappings involving complementarity relations. And in the final section we prove Theorems \ref{rob} and \ref{dr96}. We conclude the paper by  brief bibliographic comments.

\noindent{\bf Some notation:}

$K$ stands for a cone (convex, usually polyhedral) and $L$ for a subspace;

$C^{\circ}=\{y:\; \lan y,x\ran\le 1\}$ is the polar of $C$, in particular,

$K^{\circ}=\{ y:\; \lan y,x\ran\le 0,\; \forall x\in K \}$  is the polar cone of  $K$;

$L^{\perp}=\{y:\; \lan y,x\ran= 0  \}$ is the orthogonal complement to $L$;

${\rm lin}~\!K= K\cap (-K)$ is the linearity of $K$ (the maximal linear subspace contained in $K$);

${\rm ri}~\!C$ is the relative interior of a convex set $C$;

$L(C)$ is the subspace spanned by $C-C$;

${\rm aff}~\!C$ is the affine hull of $C$;

$d(x,C)$ is the distance from $x$ to $C$.

$T(C,x)=\{ h:\;\dis\liminf_{t\to+0}t^{-1}d(x+th,C)=0 \}$ is the tangent (contingent) cone to $C$ at $x\in C$;

$N(C,x)= \{y:\; \lan y,u-x\ran\le 0,\; \forall \; u\in C \}$ is the normal cone to $C$ at $x$.

\section{A brief excursion into the regularity theory of variational analysis.}

A set-valued mapping $F: X\rra Y$ is {\it (metrically) regular near} $\xyb\in\gr F$ if there is a $\kappa>0$ such that $d(x,F^{-1}(y))\le \kappa d(y,F(x))$ for all $(x,y)$ close to $\xyb$. We say that $F$ is {\it globally} metrically regular
or just {\it regular} if the inequality holds for all $(x,y)\in\R^n\times\R^n$.
The lower bound $\reg F(\xb|\yb)$ of such $\kappa$ is the {\it modulus} or {\it  rate of metric regularity} of $F$ near $\xyb$. If no such $\kappa$ exists, we set
$\reg F(\xb|\yb)=\infty$.

An equivalent property: $F$ is  {\it open at a linear rate near} $\xyb$ if there are $r>0$
and $\ep >0$ such that $B(y,rt)\cap B(\yb,\ep)\subset F(B(\xb,t))$ for all $(x,y)\in\gr F$ sufficiently close to
$\xyb$ and $t\in(0,\ep)$.  The upper  bound $\sur F(\xb|\yb)$ of such $r$ is the {\it modulus} or {\it rate of surjection} of $F$ near $\xyb$. If no such $r$ exists, we set
$\sur F(\xb|\yb)= 0$. The equality
$$
\reg F(\xb|\yb)\cdot \sur (F(\xb|\yb))=1
$$
always holds, provided we adopt the convention that $0\cdot\infty = 1$.

 There is a third
equivalent description of regularity, through the so called {\it pseudo-Lipschitz}
or {\it Aubin} property of the inverse map. We do not state it in full here, just mention that in case
when $F$ is regular near $\xyb$ and the restriction of the inverse mapping $F^{-1}$ to a neighborhood of $\xyb$ is single-valued, it is necessarily Lipschitz.

We refer to a recent survey \cite{AI15} for an account of an essential part of the metric regularity theory. Here
we just mention a few facts needed for our purpose.

\begin{definition}\label{streg}{\rm
		The set-valued mapping $F: X\rra Y$ (both $X$ and $Y$ metric) is {\it strongly
			regular} near $\xyb\in\gr F$ if it is metrically regular near $\xyb$ and there is a
		$\del >0$ such that $F(x)\cap F(u)\cap B(\yb,\del)=\emptyset$ for $x,\ u$ sufficiently
		close to $\xb$. If $F$ is globally regular and $F(x)\cap F(u)=\emptyset$ for all $x\neq u$, we say that $F$ is {\it globally} strongly regular.
	}
\end{definition}

Strong regularity is the key property behind practically all uniqueness and parametric stability  theorems for solutions of various inclusions, generalized equations and
variational inequalities, in particular. The reason is that $F$ {\it is strongly regular
	near $\xyb\in\gr F$ if and only if for some $\ep>0$ the restriction of
	$F^{-1}(y)\cap B(\xb,\ep)$ to a neighborhood of $\yb$ is single-valued and Lipschitz.}
It is through this  property that strong regularity was originally defined
by  Robinson in \cite{SMR80}.

%

A set-valued mapping $F: X\rra Y$ between linear spaces is {\it homogeneous} if its graph is a pointed cone. This means that $F(\la x)=\la F(x)$ and $F^{-1}(\la y)=\la F^{-1}(y)$.

\begin{proposition}\label{globa}
Let $X$ and $Y$ be two Banach spaces, and let $F: X\rra Y$ be a homogeneous set-valued
mapping. If $F$ is (strongly) regular near $(0,0)$, than it is globally (strongly) regular with the same rates.
\end{proposition}

\proof
By the assumption, there are $\kappa>0$ and $\ep >0$ such that
$d(x,F^{-1}(y))\le \kappa d(y,F(x))$ if $\max\{\| x\|,\| y\|  \}<\ep$
(and in case $F$ is strongly regular, $F(x)\cap F(u)=\emptyset$ for such $x$ and $u$ if $x\neq u$). Take now arbitrary $(x,y)\in X\times Y$. Set $\| m\|= \max\{\| x\|,\| y\|  \}$, and let
$\mu<\ep/m$. Then
$$
\mu d(x,F^{-1}(y))=d(\mu x,F^{-1}(\mu y))\le \kappa d(\mu y, F(\mu x))=\mu\kappa d( y,F(x))
$$
whence $d(x,F^{-1}(y))\le \kappa d(y,F(x))$. If $F$ is strongly regular near $(0,0)$, then
$F(\mu x)\cap F(\mu u)=\emptyset$ if $x\neq u$ and $\mu$ is sufficiently small, hence $F(x)\cap F(u)=\emptyset$. \endproof

The proposition allows to get   simplified  formulas for the regularity rates for
closed-graph homogeneous mappings between finite dimensional spaces. Namely

\begin{proposition}\label{cod}
Let $F: \R^n\rra \R^m$ be a homogeneous set-valued mapping with closed graph. Then
    $$
    \sur F(0,0)=\inf_{(x.y)\in\gr F}\inf\{\| x^*\|:\; (x^*,y^*)\in[T(\gr F,(x,y))]^{\circ},\; \| y^*\|=1   \}.
    $$
\end{proposition}

\section{Polyhedral sets and mappings.}
        A set $C\subset \R^n$ is {\it polyhedral}
        if it is an intersection of finitely many closed linear half spaces,
        that is if
        \begin{equation}\label{polyd}
        C=\{x\in\R^n:\; \lan y_i,x\ran\le\al_i,\ i=1,\ldots,k \}
        \end{equation}
        for some nonzero $y_i\in\R^n$ and  $\al_i\in\R$. Note that some of the inequalities may actually hold as equalities for elements of $C$, if $y_i=-y_j$ for some $i$ and$j$.

    If all $\al_i$ are equal to zero, we get a definition of a {\it polyhedral cone}.
        A set-valued mapping $\R^n\rra\R{^m}$ is {\it polyhedral} if its graph is polyhedral.
Clearly, any polyhedral set is closed and convex.

\begin{definition}\label{face}{\rm
Given a closed convex  set $C\subset \R^n$, a closed subset $F\subset C$ is a {\it face} of
$C$, if for any line segment $[a,b]$
$$
[a,b]\subset C\;\&\; x\in{\rm ri}~\![a,b]\;\Rightarrow\; [a,b]\subset F.
$$
We denote by $\cf_C$ the collection of all faces of $C$ and set
for $x\in C$, $\bar F\in\cf_C$
$$
\cf_C(x)=\{F\in\cf_C:\; x\in F\},\quad \cf_C(\bar F)=\{F\in\cf_C:\; \bar F\subset F  \}.
$$

A face $F\in\cf_C$ is {\it exposed} if there is a nonzero $y\in \R^n$ such that
$$
F=\{ x\in C:\; \lan y,x\ran = \max_{u\in C}\lan y,u\ran\}
$$
We shall say that $y$ {\it exposes F} in $C$.
}
\end{definition}

Clearly, given a representation of $C$
as in (\ref{polyd}), every face of $C$ is a collection of elements of $C$
satisfying equations $\lan y_i,x\ran=\al_i$ for some $i$.
Therefore any face is itself a polyhedral set. The biggest face is the set itself. We shall call a face of $C$ {\it proper} if is different from $C$.
The following
facts about faces are well known (see \cite{RTR}, \S\S~18,19) and geometrically obvious.

\begin{proposition}\label{polyprop} Let $C$ be a polyhedral set given by (\ref{polyd})

(a) for any $I\subset{i=1,\ldots,k}$ the set $F_I=\{x\in C:\; \lan y_i,x\ran=\al_i,\; i\in I  \}$, if nonempty,
is a face, and for any $F\in\cf_C$ there is an $I$ such that $F=F_I$.
Hence $C$ has finitely many faces and
intersection of two faces (if nonempty) is a face;

(b) \ ${\rm ri}~\!F_I=\{x\in F_I:\; \lan y_i,x\ran< 0,\; \forall\; i\not\in I  \}$.
The family $\{{\rm ri}~\!F:\; F\in\cf_C \}$ is a partition of $C$;

(c) \ if $F\in\cf_C$ and $F'\in\cf_F$, then $F'\in\cf_C$; if $F,F'\in\cf_C$ and  $F\cap {\rm ri}~\!F'\neq\emptyset$, then $F'\subset F$;


(d) \ $C\cap {\rm aff}~\!F= F$ for any $F\in\cf_C$; every face of $C$ is exposed;

(e) \ if $C$ is a cone, then any face of $C$ is a cone; ${\rm lin}~\!F= {\rm lin}~\!C$
for any $F\in \cf_C$, that is, if $C$ is a cone, ${\rm lin}~\!C$ is the minimal face of $C$ .
\end{proposition}

As $\cf_C$ is a finite set, it follows from (a) that for any $x\in C$ there is a minimal face containing $x$. We shall denote it $F_{min}^C(x)$. Given a representation of $C$ by (\ref{polyd})
and $x\in C$, let $I(x)=\{i\in \{1,\ldots,k\}:\; \lan y_i,x\ran= \al_i\}$. Then
$$F_{min}^C(x)= \{u\in C:\; \lan y_i,u\ran=\lan y_i,x\ran,\; i\in I(x) \} =F_{I(x)}.$$
In view of (b), the following statement is now obvious.

\begin{proposition}\label{minface}
Let $C$ be a polyhedral set. If $F\in\cf_C$, then  $x\in{\rm ri} ~\!F$ if and only if $F=F_{min}^C(x)$.
\end{proposition}

Elementary geometric argument allow to reveal  that the orthogonal
projection of a polyhedral set and, more generally, a linear image
of a polyhedral set is also a polyhedral set (see \cite{RTR},
Theorem 19.3).

\begin{proposition}[local tangential representation]\label{loctan}
    Let $C\subset \R^n$ be a polyhedral set and $\xb\in C$. Then there is
    an $\ep >0$ such that
    $$
    C\cap B(\xb,\ep)= \xb+ T(C,\xb)\cap(\ep B).
    $$
\end{proposition}

\proof  Let
$C=\{ x:\; \lan x_i^*,x\ran\le \al_i,\; i=1,\ldots,k\}$.
It is an easy matter to see that
\begin{equation}\label{loctan1}
T(C,\xb)=\{h:\; \lan x_i^*,h\ran\le 0,\; i\in I(\xb)  \}.
\end{equation}
It follows that $\xb+h\in C$ for an $h\in T(C,\xb)$ if $\lan x_i^*,x+h\ran<\al_i$
for all $i\not\in I(\xb)$. Clearly, the latter holds if $\| h\|<\ep$ with
sufficiently small  $\ep>0$. As $C\subset \xb +T(C,\xb)$ (which is true for any closed convex set, not necessarily polyhedral), the result follows.
\endproof

The proposition says that the local geometry of a polyhedral set or of a set which is a finite union of polyhedral sets near a certain point is fully determined by the geometry of its tangent cone at the point.
Thus,
the regularity properties of a  set-valued mapping  whose graph is a finite union of
polyhedral sets near a point of the
graph are fully determined by the corresponding properties of the graphical derivative
 of the mapping at the point (the set-valued mapping whose graph is the tangent cone to the graph - see \cite{RW}) near $(0,0)$. This often allows to work with polyhedral cones rather than with general polyhedral sets.

Needless to say that the tangent cone to a polyhedral set at any point of the set is a polyhedral cone.
A useful and equally elementary property of tangent cones to polyhedral sets
partly following from the last proposition is

\begin{proposition}\label{tandif}
Let $C\subset \R^n$ be a polyhedral set and $x\in C$. Then
$$
T(C,x) = \cone (C-x) = \cone (C-F_{min}^C(x)).
$$
Moreover, if $C$ is a cone, then $T(C,x)= C+L(F_{min}^C(x))$.
\end{proposition}

\proof The first equality in the displayed formula is an easy consequence of Proposition \ref{loctan}. To prove the second, take some $u\in C$ and $v\in F_{min}^C(x)$. Replacing, if necessary, $u$ by
$x+\la(u-x)$ and $v$ by $x+\la(v-x)$ (this does not change the direction $u-v$), we may assume that $w= x+ (x-v)\in C$ (as $x\in\ {\rm ri}~\! F_{min}^C(x)$). Then
$z=(1/2)(u+w)\in C$ and therefore $z-x\in C-x$. But
$2(z-x) = u-v$, whence $u-v\in\cone (C-x)$. This proves that $\cone (C-F_{min}^C(x))\subset\cone(C-x)$. The opposite inclusion is trivial.
The cone sign can be dropped in case $C$ is a cone as $C-F$ is a cone for any $F\in\cf_K$ and
the simple chain $K-F= (K+F)-F= K+L(F)$ proves the last equality.
\endproof

In turn, this immediately implies
\begin{proposition}\label{normface}
	Let $C$ be a polyhedral set and $F\in\cf_C$. If $x,x'\in{\rm ri} ~\!F$, then
	$I(x)= I(x')$, $T(C,x)=T(C,x')$ and $N(C,x)=N(C,x')$.
\end{proposition}

This justifies introducing notation $I(F)$, $T(C,F)$ and $N(C,F)$ for the common index sets, and tangent and
normal cones to $C$ at elements of the relative interior of a face. It is natural to call
them the {\it tangent} and the {\it normal cones to $C$ at $F$}.

Finally, it follows from the standard rules of convex analysis (see e.g. \cite{DR}, Theorem 2E.3) that
for $x\not\in \intr C$ the normal cone $N(C,x)$ is the convex cone generated by
$\{y_i:\; i\in I(x) \}$.  It follows from the proposition that, given a face $F$ of $C$,  cones
\begin{equation}\label{srule}
N(C,F) = \{\sum \la_iy_i:\; i\in I(F),\; \la_i\ge 0  \}.		
\end{equation}

\section{Face complementarity}
\begin{definition}\label{dual}{\rm
		Let $K$ and $H$ be two polyhedral cones in $\R^n$. We say that $K$ and $H$ are {\it face complementary}  if there is a
		one-to-one correspondence $\Lambda$ between $\cf_K$ and $\cf_H$, such that for any $F,\ F'\in\cf_K$.
		
		
		(i) \ \quad $\dim L(F)+\dim L(\Lambda(F))= n$, 
		
		(ii) \quad $F\subset F'\;\Rightarrow\; \Lambda (F')\subset \Lambda(F)$.
		
		\noindent  We also say that the complementarity correspondence is {\it non-singular}  if
		\begin{equation}\label{lin1}
		    L(F)\cap L(\Lambda(F))= \{0\};\quad          L(F)+ L(\Lambda(F))= \R^n
		\end{equation}		
	}
\end{definition}

In what follows we write $\Lambda_K$ and $\Lambda_H$ to specify the direction in which we apply $\Lambda$, namely $\Lambda_K: \cf_K\to\cf_H$ and $\Lambda_H: \cf_H\to\cf_K$.

The well known  example of face complementary polyhedral cones is $K$ and $K^{\circ}$ (see e.g. \cite{SMR92}, Proposition 2.1). We shall give a short proof of this fact later in this section. Of course, $\Lambda$ may  not be uniquely defined -- there may be other face complementarity correspondences between $K$ and $H$. But we assume throughout that, given $K$ and $H$,  a  certain $\Lambda$ is fixed.
The collection of face complementary pairs of polyhedral cones is rather rich.
As a trivial example,
just mention that  two polyhedral  cones in $\R^3$   are face complementary if and only if they have the same number of extreme directions.

\begin{proposition}\label{invnorm}
	Let $K$ and $H$ be a pair of face complementary polyhedral cones, and let
	$\Lambda$ be the corresponding face complementarity correspondence. Then $\Lambda_H(F_{min}^H(y))$ is the maximal among faces $F\in\cf_K$ such that $y\in \Lambda(F)$. Likewise $\Lambda_K(F_{min}^K(x))$
	is the maximal element among $G\subset \cf_H$ such that $x\in\Lambda_H(G)$.
	Moreover the relations
	\begin{equation}\label{4.2}
	y\in \Lambda_K(F_{min}^K(x))\quad {\rm and}\quad x\in\Lambda_H(F_{min}^H)(y)
	\end{equation}
	are equivalent.
\end{proposition}
\proof Suppose $x\in F$ and $y\in \Lambda_K(F)\subset \Lambda_K(F_{min}^K(x))$. Then
$$
x\in F=\Lambda_H(\Lambda_K(F))=\Lambda_H(G)\subset \Lambda_H(F_{min}^H(y))).
$$		
This implies that $F\subset \Lambda^{-1}(F_{min}^H)(y)$ and the left inclusion in (\ref{4.2}) implies the right one. Reversing the roles of $x$ and $y$ and $K$ and $H$,
we complete the proof.\endproof

In view of this result, it is natural to set
$$
F_{max}^K(y) = \Lambda_H(F_{min}^H(y));\qquad F_{max}^H(x)= \Lambda_K(F_{min}^K(x)).
$$
We can view $F_{max}^K$ and $F_{max}^H$ as set valued mappings from $\R^n$ into itself
with domains $K$ and $H$ respectively. The proposition shows that {\it the mappings
	are mutually inverse:}
\vskip 1mm

\centerline{$\big(F_{max}^K\big)^{-1}= F_{max}^H$.}
\vskip 2mm

There are two simple transformations that preserve complementarity of polyhedral cones
and play a substantial role in further discussions.
Let $\bar F$ be a proper face of $X$ and $L=\bar F-\bar F$ the subspace spanned by $\bar F$. Set $K_{\bar F}=T(C,\bar F)$, $H_{\bar F}=\Lambda_K(\bar F)$.
By Proposition \ref{tandif} $K_{\bar F}= \cone(K-\bar F)= K+L$.
 Clearly any face of $K_{\bar F}$ has the form
$F+L$, with $F\in\cf_K(\bar F)$ and any such set is a face of $K_{\bar F}$. Define $\Lambda_{K_{\bar F}}$ by $\Lambda_{K_{\bar F}}(F+L)= \Lambda(F)$. Then
$\Lambda_{K_{\bar F}}$ is a one-to-one correspondence between $\cf_{K_{\bar F}}$ and $\cf_{H_{\bar F}}$.
It is immediate from the definitions that $\Lambda_{K_{\bar F}}(G)\subset\Lambda_{K_{\bar F}}(G')$
if $G, G'\in\cf_{K_{\bar F}}$ and $G'\subset G$. Finally, for any $F\in\cf_K({\bar F})$
the dimensions of $F$ and $F+L$ obviously coincide as $\bar F\subset F$.
Thus we have proved

\begin{proposition}\label{tantran}	
	The cones $K_{\bar F}$ and $H_{\bar F}$ are face complementary and  $\Lambda_{\bar F}$ is a complementarity correspondence between $K_{\bar F}$ and $H_{\bar F}$.
\end{proposition}

We shall say that $(K_{\bar F},H_{\bar F},\Lambda_{\bar F})$ is the {\it tangential extension}
of $(K,\Lambda_K)$ along  $\bar F$. Tangential extensions of $H$ are defined in a similar
symmetric way.

The proposition opens gates for extension of the complementarity concept to arbitrary polyhedral sets. Take such a $C$, and let $\bar F\in\cf_C$. Consider the tangent cone
$T(C,\bar F)$. Then $L={\rm lin}~\! T(C,\bar F)= \cone (\bar F-\bar F)$ and any face $G$ of $T(C,\bar F)$ has the form $F+L$, where $F\in \cf_C(\bar F)$. (The latter because
for any $\xb\in\ri\bar F$ any $x\in C$ sufficiently close to $\xb$ must belong
to a face containing points of $\bar F$.)

\begin{definition}\label{dual1}{\rm
Let $C\in\R^n$ be a polyhedral set. The set valued mapping
$\Lambda_C: \cf_C\rra \R^n$
is a {\it complementarity relation on $C$} if for any
$\bar F\in\cf_C$

\vskip 2mm

(i) \  \  $\Lambda_C(\bar F$)  is a polyhedral cone;

\vskip 2mm

(ii) \ $K_{\bar F}=T(C,\bar F)$ and $H_{\bar F}= \Lambda_C(\bar F)$ are face complementary cones
and relations
$$
(\Lambda_{\bar F})_{K_{\bar F}}(\cone(F-\bar F))=\Lambda_C(F),\quad (\Lambda_{\bar F})_{H_{\bar F}}(G)= \cone(F-\bar F)
$$
for $F\in\cf_C(\bar F),\; G\in\Lambda_C(F)$ define a complementarity correspondence between $K_{\bar F}$ and $H_{\bar F}$.
We shall call the complementarity relation $\Lambda_C$ {\it non-singular} if
$L(F)+L(\Lambda_C(F))=\R^n$ (and consequently $L(F)\cap  L(\Lambda_C(F))=\{0\}$) for all $F\in\cf_C$.
}

\end{definition}

Returning back to the case of two complementarity cones, we can now say that
$\Lambda_K$ is a complementarity relation on $K$ and $\Lambda_H$ is a complementarity
relation on $H$. Likewise the triple $(K_{\bar F},H_{\bar F},\Lambda_{\bar F})$ in the definition can be defined as
the {\it tangential extension of $(C,\Lambda_C)$ along $\bar F$}.

Now we can describe the second complementarity preserving transformation, this time
already for any polyhedral set $C$. Let  $\Lambda_C$ be a non-singular complementarity relation on $C$, and let $\bar F\in\cf_C$. Denote for brevity $M=L(\Lambda_C(\bar F))$, the subspace spanned by $\Lambda_C(\bar F)$.  If $\Lambda_C$ is non-singular, $M$ nd $L=L(\bar F)$
are complementary subspaces of $\R^n$ and we can consider the projection $\pi_{ML}$
onto $M$ parallel to $L$.

Set $K'_{\bar F} =\cone(\pi_{ML}(C))=\pi_{ML}(T(C,\bar F))=M\cap T(C,\bar F)$, \
$H'_{\bar F}=\Lambda_C(\bar F)$ (note that $H'_{\bar F}$ coincides with $H_{\bar F}$
introduced in Definition \ref{dual1} above). It is an easy matter to see that
any face of $K'_{\bar F}$ is the image under $\pi_{ML}$ of a face of $T(C,\bar F)$
and any face of $K'_{\bar F}$ can be obtained this way. Thus $G\in\cf_{K'_{\bar F}}$
if and only if there is an $F\in\cf_C(\bar F)$ such that $G=\cone(\pi_{ML}(F-\bar F))$.
Define a one-to-one correspondence between faces of $K'_{\bar F}$ and $H'_{\bar F}$
by setting for $F\in\cf_C(\bar F)$ and $G\in\Lambda_C(\bar F)$
$$
(\Lambda_{\bar F})_{K'_{\bar F}}(\cone(\pi_{ML}(F-\bar F))=\Lambda_C(F),\quad
(\Lambda_{\bar F})_{H'_{\bar F}}= \pi_{ML}(cone(F-\bar F).
$$

\begin{proposition}\label{proj}
Let $C$ be a polyhedral set in $\R^n$ and $\Lambda_C$ a non-singular complementarity relation on $C$. Then for any $\bar F\in\cf_C$ the cones $K'_{\bar F}$ and
$H'_{\bar F}$ are face complementary in $M$ with $\Lambda_{\bar F}$ a face complementarity correspondence between $K'_{\bar F}$ and
$H'_{\bar F}$.
\end{proposition}

\proof Immediate consequence of Proposition \ref{tantran}.

We shall call the triple $K'_{\bar F},H'_{\bar F},\Lambda'_{\bar F}$ the
{\it factorization of $(C,\Lambda_C)$ along $\bar F$.}

\begin{proposition}\label{nonun}
	Let $C\subset \R^n$ be a polyhedral set, and let $\Lambda$ be a complementarity relation on $C$. If
	$A$ is a linear isomorphism $\R^n\to\R^n$, then $\Lambda_A(A(F))=\Lambda(F)$ is a complementarity relation on $A(C)$.
\end{proposition}

The proposition says that the complementarity property is stable with respect to linear isomorphisms.
We omit the proof of the proposition which is fairly straightforward. The only remark that probably should be made is that  the values of $A$ outside of $C$ are inessential for the result and it is enough to assume that $A$ is injective on $L(C)$.

To conclude the section, we consider the most important example of face complementarity
that will be our primal interest when we come back to variational inequalities.

\begin{proposition}\label{comnorm}
	Let $C\subset \R^n$ be a polyhedral set. Then $\Lambda_C: F\to \N(C,F)$
is a non-singular complementarity relation on $C$. In particular, if  	
	$K$ is a
	polyhedral cone in $\R^n$, then
	$K$ and $K^{\circ}$ are face complementary with $\Lambda:
	F\leftrightarrow N(K,F)$ being a non-singular face complementarity correspondence.
\end{proposition}
\proof
Note that the first statement is an easy consequence of the second
in view of Proposition \ref{loctan}. To prove the the second statement, we
observe first that $N(K,F)$ is a face of $K^{\circ}$. Indeed,
it is clear that $N(K,F)\subset K^{\circ}$. On the other hand, given
an $F\in\cf_K$ and  $x\in {\rm ri}~\!F$, we have $\lan y,x\ran=0$
for  $y\in N(K,F)$ and $\lan v,x\ran \le 0$ for $v\in K^{\circ}$.
Thus, if $y=(1/2)(y_1+y_2)$, $y_i\in K^{\circ}$, we have $\lan y_i,x\ran=0$,
that is $y_i\in N(K,x)$.

It follows from (\ref{srule}) and Proposition \ref{polyprop}(a) that
$N(K,F_1)\neq N(K,F_2)$ if $F_1\neq F_2$.  For the same reason
$F_1\subset F_2$ implies $N(K,F_2)\subset N(K,F_1)$. Finally,
if  $K=\{x:\; \lan y_i,x\ran\le 0,\; i=1,\ldots,k\}$ and
$F=F_I$ for some $I\subset\{1,\ldots,k\}$, then $L(N(K,F))$ is the
linear hull of $\{ y_i,\; i\in I\}$ which means that $L(F)$ and
$L(N(F))$ are complementary subspaces of $\R^n$, hence the sum of
their dimensions and the dimension of their sum is $n$.\endproof

It is an easy matter to see that $N(K,\cdot)$ and $N(K^{\circ},\cdot)$ are mutually inverse
mapping that is for
\begin{equation}\label{4.1}
G=N(K,F) \;  \Leftrightarrow\; F= N(K^{\circ},G).
\end{equation}

\section{The complementarity map}

Let $C\subset \R^n$ be a polyhedral set and $\Lambda_C$ a complementarity relation on $C$.  The {\it complementarity mapping associated with $(C,\Lambda_C)$}  is
\begin{equation}\label{maps}
\Phi_C(x) = x+\Lambda_C(F_{min}^C(x)). 
\end{equation}

In case we have two face complementary cones $K$ and $H$ with face complementarity correspondence $\Lambda: \cf_K\leftrightarrow\cf_H$, then $\Phi_K(x) = x+F_{max}^H(x)$
and we can consider the symmetric mapping $\Phi_H(y)= y+ F_{max}^K(y)$.
By Proposition \ref{invnorm}
\begin{equation}\label{sym}
x\in F_{max}^K(y) \Leftrightarrow\; y\in F_{max}^H(x)\;\Leftrightarrow\;
z=x+y\in\Phi_H(x)\cap\Phi_H(y).
\end{equation}
As a consequence of this observation, we get that the ranges of $\Phi_K$ and $\Phi_H$
coincide and, moreover, the following holds.
\begin{proposition}\label{reg}
	If $\Phi_K$ is (strongly) regular near $(\xb,\zb)\in\gr\Phi_K$,
	then $\Phi_H$ is (strongly) regular near $(\yb,\zb)$, where $\yb=\zb-\xb$.
\end{proposition}

\proof Let $(x,z)\in\gr \Phi_H$ be sufficiently close to $(\xb,\zb)$. As $\Phi_K$
is regular near $(\xb,\zb)$, for any $w$ sufficiently close to $z$ there is a $u$
such that $w\in\Phi_K(u)$ and $\| u-x\|\le r^{-1}\| z-w\|$, where $r$ is the
rate of surjection of $\Phi$ near $(\xb,\zb)$. We have $z=x+y, \ w=u+v$ for some
$y\in \Lambda_K(x)$ and $v\in \Lambda_K(u)$ and therefore
$$
\| y-v\|\le \| z-w\|+ \| x-u\|\le k\| z-w\|
$$
with $k$ not depending on $z,w,y,v$. But by (\ref{sym}) $z\in \Phi_H(y)$ and $w\in \Phi_H(v)$.

Now if $\Phi_K$ is strongly regular near $(\xb,\zb)$, then there is a Lipschitz mapping
$x(z)$ from a neighborhood $W$ of $\zb$ into a neighborhood $U$ of $\xb$ such that
$\{(x,z):\; x=x(z),\; z\in W\}=(U\times W)\cap \gr\Phi_K$. Set
$y(z) = z-x(z)$. This is a Lipschitz mapping and, obviously, there is no other solutions to $z\in\Phi_H(y)$ in a small neighborhood of $\yb$ close to $\zb$.
\endproof

Both $\Phi_K$ and $\Phi_H
$ are obviously homogeneous mappings.
As a consequence we get, thanks to Proposition \ref{globa},
\begin{proposition}\label{glob}
	If $\Phi_K$ is regular near $(0,0)$,
	then it is globally regular (in particular, regular near every point of its graph) with the same rates.
\end{proposition}
The proposition allows to speak just about regularity of $\Phi_K$ and $\Phi_H$ without specifying
whether it is global or local near $(0,0)$ (but of course regularity near another
point of the graph may not imply global regularity).
This make the complementarity maps with a pair of face complementary cones especially convenient to work with.

Let us return to the general model (\ref{maps}). Fix a  $\bar F\in\cf_C$, and let
$(K_{\bar F},H_{\bar F},\Lambda_{\bar F})$ be the tangential extension of $(C,\Lambda_C)$ along  $\bar F$.
Set $\Phi_{\bar F}(x) = x+ \Lambda_C(F_{min}^{K_{\bar F}}(x))$.

\begin{proposition}\label{tantran1}
The set-valued mapping $\Phi_{\bar F}$ is (strongly) regular near $(0,\zb-\xb)$, provided  $\Phi_C$
is (strongly) regular near $(\xb,\zb)$ and vice versa. In particular, if $\Phi_C$
is (strongly) regular near $(\xb,\xb)$, then $\Phi_{\bar F}$ is globally (strongly)
regular.
\end{proposition}
\proof
 As we have seen in the proof of Proposition \ref{tantran}
for $(x,z)$ sufficiently close to $(0,\zb)$ the relation
$z-\xb\in \Phi_{\bar F}(x)$ is equivalent to
$z\in \Phi_C(x+\xb)$ This implies the first statement. The second follows from Proposition \ref{globa}. \endproof

Below are  the main results of the section.

\begin{theorem}\label{lin}
	Let $C$ be a polyhedral set, let $\Lambda_C$ be a complementarity relation on $C$, and let $\Phi_C(x) = x+\Lambda_C(F_{min}^C(x))$
	If  $\Phi_C$ is regular, then for any $F\in\cf_C$ the subspaces
	$L(F)$ and $L(\Lambda(F))$ are complementary subspaces of $\R^n$, that is
$\Lambda_C$ is a non-singular complementarity relation.
\end{theorem}

\proof
It is sufficient to consider the case of two face complementary cones. The general case
reduces to that one with the help of Propositions \ref{loctan}, and \ref{tantran1} if we consider $K=T(C,F)$ and $H=\Lambda_C(F)$.

Recall that for any $F\in\cf_K$ the sum of dimensions of $L(F)$ and $L(\Lambda(F))$
is $n$.
So
assuming that the statement is wrong, we  find an $F\in\cf_K$,  $x_i\in L(F)$ and
$y_i\in \Lambda_K(F)$, $i=1,2$ such that $0\neq x_1-x_2=y_2-y_1$. Adding, if necessary,
to both $x_i$ some $x\in {\rm ri}~\! F$ and to both $y_i$ some $y\in {\rm ri}~\! \Lambda_K(F)$, we may guarantee that $x_i\in {\rm ri}~\! F$ and $y_i\in {\rm ri}~\! \Lambda_K(F)$.
Then $F=F_{min}^K(x_i)$.

Let $(K_F,H_F,\Lambda_F)$ be the tangential extension of $(K,\Lambda_K)$
along $F$.  Set as above
$\Phi_F(x)= x+F_{max}^{H_F}(x)$.
By Proposition \ref{tantran1} $\Phi_F$ is regular, so there is an $r>0$ such that
$B(z,r\ep)\subset \Phi_F(B(x,\ep))$
for any $x\in K_F$, any $z\in \Phi_F(x)$ and sufficiently small $\ep>0$.
Set $z= x_1+y_1=x_2+y_2$. Then for any $z_i'\in B(z,\ep)$ there is a $x_i'\in K_F$
such that $z_i'\in\Phi_F(x_i')$ and $\| x_i'-x_i\|\le r^{-1}\| z_i'-z\|$.
Set further $y_i'=z_i'-x_i'\in H_F$. We have
$\| y_i'-y_i\|\le (1+r^{-1})\| z_i'-z\| \le (1+r^{-1})\ep$. Let $\ep$ be so small that
$y_i'\in{\rm ri}~\!H_F$. Then necessarily $x_i'\in L(F)$ as $L(F)$ is the smallest
face of $K_F$. This means that $\Phi_K(B(x_i,\ep))\subset F+ H_F$ and therefore
$$
B(0,2r\ep)\subset \Phi_F(B(x_1,\ep))-\Phi_1(B(x_2,\ep))\subset L(F)+L(H_F)=
L(F)+L(\Lambda_K(F)).
$$
As the sum of dimensions of $L(F)$ and $L(\Lambda(F))$ is $n$ and the subspaces
have a common nonzero vector $x_1-x_2$,  the dimension of their sum is strictly smaller than $n$ in contradiction with the displayed inclusion above. This completes the proof.
\endproof

Combining this with Propositions \ref{proj} and \ref{tantran1} we get

\begin{corollary}\label{proj1} If $\Phi_C$ is  regular, then
for any $F\in\cf_C$ the complementarity mapping $\Phi_F'$ associated with the
factorization $(K_F',H_F',\Lambda_F')$ of $(C,\Lambda)$
along $F$  is also regular in $M$.
\end{corollary}

\proof
We have $K_F'= K_F\cap M$, $H_F'=H_F$, where $(K_F,H_F,\Lambda_F)$
is the tangential extension of $C,\Lambda_C$ along $F$.
So if $z\in \Phi_F'(x)$, that is
$z=x+y$ with both  $x$ and $y$ in $M$ and $z'\in M$, then by regularity of
$\Phi_F$ (see Proposition \ref{tantran1}) there are $x'\in K_F$ and $y'\in H_F\subset M$
such that $z'= x'+y'$ and $\|x'-x\|\le \kappa\| z'-z\|$ with $\kappa$ being the modulus of regularity of $\Phi_F$. Then $x'\in M$ as both $z'$ and $y'$ are in $M$, that is
$x'=\pi_{ML}(x)$ and $z'\in \Phi_F'(x')$.\endproof

\begin{theorem}\label{sep}
Let $K$ and $H$ be a pair face complementary cones. If $\Phi_K$ is regular, then
$$
K\cap H =\{0\}.
$$	
\end{theorem}

\proof The proposition trivially holds if $n=1$. Indeed, if $K=\{0\}$, then $H$ must be a one-face  cone of dimension $1$, that is $H=\R$. If $K$ is a half line
then $H$ contains two faces, that is $H$ is also a half line. By regularity, $K$ and $H$ cannot coincide, hence zero is the only common point for $K$ and $H$.

Suppose now that the theorem is true for all dimensions up to $m-1$ with some $m\ge 2$
and consider the case $n=m$. Assume that there is a $0\neq e\in H\cap K$.
If $\dim K=k<n$, then by Theorem \ref{lin}, $\Lambda(K)$, the minimal face of $\cf_H$, must have dimension $n-k$. This may be only if it is a subspace, that is
$\Lambda_K(K)={\rm lin}~\! H$.

Consider factorization of $(H,\Lambda_H)$ along $\bar G={\rm lin}~\! H$.
Then $K_{\bar G}=K$, $H'_{\bar G}=\pi_{ML}(H)$, where $M=L(K)$,
$L={\rm lin}~\! H$ and  the complementarity
correspondence   $\Lambda_{\bar G}': \cf_{K}\leftrightarrow \cf_{H'}$ is defined by $(\Lambda_{\bar G}')_{K'}(F)=\pi_{ML}(\Lambda_C(F))$.  Applying consecutively Propositions  \ref{reg}
and \ref{proj}, we conclude that $\Phi_{K'}'=x+(\Lambda_{\bar G}')_{K'}(F_{min}^{K'}(x))$ is a regular
mapping in $M$. But $\dim M<m$, so by the induction hypothesis $K\cap H'=\{0\}$.
As $ {\rm lin}~\! H$ is a complementary subspace to $L(K)$, it follows that
$K\cap (H_{\bar G}'+{\rm lin}~\! H)=\{0\}$ and therefore that $K\cap H=\{0\}$.

So we have to consider the case $\dim K=m$. Assume that there is a $0\neq e\in K\cap H$.
The inclusion $H\subset K$ would imply that $K+H=K\neq \R^m$ contrary to the assumption that $\Phi$ is regular. This means that $H$ contains an element not belonging to $K$.
 This cannot be $-e$ for in this case ${\rm lin}~\! H\neq \{0\}$ and therefore
 $\dim K=m-\dim({\rm lin}~\! H)<m$. It follows (by convexity of $H$)  that $K\cap H$ must contain a nonzero
 element not belonging to $\intr K$ and we may assume that  $e$ is such an element.

 Let $F$ be a proper face of $K$ containing $e$. The inclusion $e\in\Lambda_K(F)$
 is impossible as $F$ and $\Lambda_K(F)$ lie in complementary subspaces of $\R^n$ by Theorem \ref{lin}. Thus $e\not\in \Lambda_K(F)$. Set $G=\Lambda_K(F)$ and  consider the
 factorization $(H_G',K_G',\Lambda_G$ of $(H,\Lambda_H$ along $G$, that is
  $K_G'=F$, $M=L(F)$ is the subspace spanned by $F$ and $L=L(G)$,   $\pi_{ML}$ is the projection onto $M$ parallel to $L$,
 $H_G'=\pi_{ML} (H)$ and $\Lambda_G'(F')=\pi_{ML}(\Lambda_K(F'))$ for
 $F'\in\cf_K(F)$.
 As $e\in M$, it follows that $e=\pi_{ML}(e)in H_G'$. Thus again we get $e\in K_G'\cap H_G'$ which, in view of Corollary \ref{proj1}
 contradicts to the induction hypothesis. This completes the proof of the theorem.
 \endproof

\section{Complementarity mapping: main results}

Everywhere in this section $C$ is a polyhedral set, $\Lambda_C$ is a complementarity relation on $C$ and $\Phi_C(x)=x +\Lambda_C(F_{min}^C(x))$ is the corresponding
complementarity mapping.

\begin{theorem}[regularity implies strong regularity]\label{main}
Let $C$ be a polyhedral set and $\Lambda$ is complementarity relation on $C$.	
 If $\Phi_C$ is regular than the
inverse mapping is single-valued and Lipschitz on $\R^n$. Thus,  regularity
of $\Phi$ implies global strong regularity.
\end{theorem}

\proof
 We only need to show that $\Phi_C^{-1}$ is single-valued: the Lipschitz property
will then automatically follow from regularity. To begin with we note that
the equality $u=x+y$ for some $x,u\in C$ and $y\in \Lambda (F_{min}^C(x))$
can be valid only if $x=u$ and $y=0$. Indeed, as $\Phi_C$ is regular,
for any $x\in C$ the complementarity mapping $\Phi_{ F}$ associated with
then tangential extension $(K_{F},H_{F},\Lambda_{F})$ of $(C,\Lambda_C)$
along $ F= F_{min}^C(x)$is also regular
(Proposition \ref{tantran1}). By Theorem \ref{sep} $K_{F}\cap H_{F}=\{ 0\}$.
If we now assume that $u=x+y$ with $x,u,y$ as above, then $u\not\in F=F_{min}^C(x)$
for otherwise $u-x=y$ and as $L(F)$ and $L(\Lambda_C(F))$ are complementary subspaces
by Theorem \ref{lin}, we get $x-u=y=0$.
But $C-x\subset T(C,\xb)=K_{F}$ and therefore for no $u\in C$ other than $x$
the difference $u-x$ may belong to $H_F=\Lambda_C(F)$.

Note that in case we have a couple of face complementary cones $K$ and $H$
with face correspondence $\Lambda$, and $\Phi_K$ is regular, then so is $\Phi_H$
(Proposition \ref{reg}) and, likewise, the equality $y=v+u$ for $y,v\in H$ and
$u\in \Lambda^{-1}(F_{min}^H(v))$ may hold only if $y=v$ and $u=0$.

 Assume now that $x+y=u+v$ for some $x,u\in C,\; x\neq u$, $y\in \Lambda (F_{min}^C(x))$
and $v\in \Lambda (F_{min}^C(u))$. We first show that this is impossible if in addition
$x\in F_{\min}^C(u)$. Indeed, in this case $F=F_{\min}^C(x)\subset F_{\min}^C(u)$ and therefore
$v\in\Lambda(F_{min}^C(x))$. Set $F=F_{min}^C(x)$, and
let $(K_F',H_F',\Lambda_F')$ be the factorization of
$(C,\Lambda)$ along $F$. Then $y,v\in H_F'$,  $u'=\pi_{ML}(u)$ (where as usual $\pi_{ML}$ is
the projection onto $M=L(\Lambda_C(F))$ along $F$) and $y= u'+v$. But by what we have just proved  this can be true only if $y=v$.

Thus we have to assume that neither of $x,\ u$ belongs to the minimal face of the other.
Let $\kappa$ be the modulus of metric regularity of $\Phi$ or any bigger number.
Choose $\ep>0$ so small that the ball of radius $(1+\kappa)\ep$ around $x$ does not meet any face
$F\in\cf_C$  such that  $({\rm ri}~\!F_{min}^C(x))\cap F=\emptyset$ .
This means that $x\in F_{min}^C(w)$ whenever $w\in C$ and $\| w- x\|\le (1+\kappa)\ep$.
Let further $m$ be an integer
big enough to guarantee that $\del= m^{-1}\| y\|<\ep$. Regularity of $\Phi_C$ allows to construct recoursively
a finite sequence of pairs $(u_k,z_k), \ k=0,1,\ldots,m$ such that
$$
(u_0,z_0)=(u,z),\quad z_k\in F_{max}^K(u_k),\quad u_k+z_k= x+(1-m^{-1}k)y,\quad \|u_k-u_{k-1}\|\le \kappa\del.
$$
Then $u_m+z_m=  x$. As follows from the result obtained at the beginning of the proof,
this can happen only if $u_m=x$. This in turn means, if $u_0\neq x$,  that for a certain $k$ we have $u_k\neq x,\; \| u_k-x\|\le \kappa\del< \kappa\ep$. By the choice of $\ep$ this implies that
$x\in F_{min}^C(u_k)$. But in this case the result obtained at the preceding paragraph
excludes the possibility of the equality $u_k+z_k= x+(1-m^{-1}k)y$
unless $u_k=x$. So we again get a contradiction that completes the proof.\endproof

Our nex goal is to characterize regularity of $\Phi_C$ in verifiable terms.

\begin{definition}{\rm
We say that $(C,\Lambda_C)$ satisfies the {\it face separation condition} if
 $\Lambda_C$ is non-singular and for any $F,F'\in\cf_C$, $F\subset F'$,
 $\dim F'=\dim F+1$,  and the hyperplane $L(F)+L(\Lambda_C(F'))$
 properly separates $L+\Lambda_C(F)$ and $F'$.}
\end{definition}
Proper separation means that none of the separated sets lies in the separating hyperplane.

\begin{theorem}[characterization of regularity]\label{chareg}
	The following properties are equivalent:
	
	(a)  $\Phi_C$ is regular;
	
	(b)  $\Lambda_C$ is non-singular and $\Phi_C$ is locally open
	
	(c) $\Lambda_C$ is non-singular and $\Phi_C$ is onto (that is the range of $\Phi$ is the whole of $\R^n$);
	
	(d) the face separation condition holds for $(C,\Lambda_C)$.
\end{theorem}

\proof (a) $\Rightarrow$ (b). Non-singularity follows from Theorem \ref{lin} and local openness is a trivial consequence of global regularity. The implication (b) $\Rightarrow$ (c) is also elementary. Indeed, local openness implies that
${\rm Im}~\!\Phi_C= \Phi_C(C)$ is an open set. On the other hand this set is the union of finitely many polyhedral sets $F+\Lambda_C(F),\; F\in\cf_C$, hence it is closed.
Being both open and closed and nonempty, it coincides with all of $\R^n$.

(c) $\Rightarrow$ (a). By non-singularity, for any $F\in\cf_C$ the mapping $(x,y)\to x+y$ from
$L(F)\times L(\Lambda_C(F))$ into $\R^n$ is a linear homeomorphism.
This means that there are positive $r(F),\; F\in\cf_C$ such that
$B(z,r(F)t)\subset B(x,t)+B(y,t)$ whenever $z=x+y$. For
$x\in F, \ z\in\Phi_C(x)$ this leads to

\vskip 2mm

\centerline{$B(z,r(F)t)\cap \Phi_C(F)\subset B(x,t)\cap F+\Lambda_C(F)\subset
	\Phi_C(B(x,t)) $ }.
\vskip 2mm

\noindent Setting $r=\min_{F\in\cf_C}r(F)$, we get

\vskip 2mm

\centerline{$B(z,rt)\subset\displaystyle{\bigcup_F} (B(z,r(F)t)\cap \Phi_C(F))\subset
	\Phi_C(B(x,t)) $ }

\vskip 2mm

\noindent which is precisely regularity of $\Phi_C$.

(b) $\Rightarrow$ (d). As $\Lambda$ is non-singular, all sets $Q(F)=F+\Lambda(F)$ have
dimension $n$, that is $\intr Q(F)\neq \emptyset$ for any $F\in\cf_C$. Take a certain $F\in\cf_C$, and let  $G\in\cf_{\Lambda_C(F)}$ be a proper face of $\Lambda_C(F)$ with
$\dim G=\dim\Lambda_C(F)-1$. Then $G=\Lambda_C(F')$ for some $F'\in\cf_C$ with
$F\subset F'$ and $\dim F'= n-\dim G= \dim F+1$. We have to show that $E={\rm aff}(F+G)$
(which is a hyperplane of dimension $n-1$) separates $Q(F)$ and $F'$.

Take an $\xb\in\ri F$ and a $\yb\in\ri G$. Then
$\zb=\xb+\yb\in Q(F)\cap Q(F')$ and $\zb\in\ri (F+G)$. We claim that
$x\in F\cup F'$ if $(x,z)\in\gr \Phi_C$ is sufficiently close to $(\xb,\zb)$.
Indeed, set $y=z-x$. Then $y\in \Lambda_C(F_{min}^C(x))$. As
$\xb\in\ri F$, any $x$ close to $\xb$ is either in $\ri F$ or in some face
containing $F$. In the last ase by the property (ii) of Definition \ref{dual} any $y\in \Lambda(F_{min}^C(x))$ belongs to $\Lambda_C(F)$. On the other hand, as $\yb\in \ri\Lambda_C(F')$, any such $y$ close to
$\yb$ must be either in $\ri \Lambda_C(F')$ in which case $x$ must be in $F'$
or belongs to a face of $\Lambda_C(F)$ containing $\Lambda_C(F')$ in which case
(again by the property (ii) of Definition \ref{dual})  $x$
belongs to some face of $C$ contained in $F'$ but different from the latter. But as we have seen, $F$ is the only such face to which $x$ may belong.

As $\Phi$ is is locally open, it follows that $Q(F)\cup Q(F')$ covers a neighborhood of $\zb$.
But $E$ is common boundary of $Q(F)$ and $Q(F')$. (Indeed $\zb$ is a boundary point of
$Q(F)$ because  $\yb$ is a boundary point of
$\Lambda_C(F)$ and, likewise, $\zb$ is a boundary point of $Q(F')$ because $\xb$
belongs to the boundary of $F'$.) Thus $Q(F)$ and $Q(F')$ must lie on different sides of $E$. But this may happen only if $E$ separates $F'$ and $F+\Lambda_C(F)$.

(d)  $\Rightarrow$ (b). We need the following simple Lemma to proceed.

\begin{lemma}\label{lem}
Let $Q_1,\ldots,Q_k$ be  polyhedral sets. Assume that $Q=\cup Q_i\neq\R^n$. Then there
are a set $Q_j$, an $\xb$ in the boundary of $Q_j$ and an $\ep >0$ such that
$B(\xb,\ep)\cap Q_j= B(\xb,\ep)\cap Q$.
\end{lemma}

\proof
The lemma is obvious if $n=1$. Suppose it is true for dimensions up to
$m-1$ for some $m\ge 2$, and let $Q_j\subset\R^m$. We can assume without loss of generality that all $Q_i$ are cones. Otherwise we can take a boundary point
of $Q$ and consider the tangent cones to $Q_i$ at $x$ instead of $Q_i$ themselves. Clearly,
zero must be a boundary point of the union of the tangent cones (immediate from Proposition \ref{loctan}).

As $Q$ is closed and $Q\neq \R^m$, there is a nonzero
boundary point $x\in Q$. Let $E$ be the subspace of elements of $\R^m$
orthogonal to $x$. Set $\hat Q_j= T(Q_j\cap(x+E))$ and $\hat Q=\cup\hat Q_i$. Clearly, in a neighborhood of $x$ all $Q_i$ coincide with the cone generated by $\hat Q_i$.
Likewise obvious is that $x$ is a boundary point of $\hat Q$ in $x+E$.

Set
$j\in J(x)=\{j:\; x\in Q_j\}$. Then   $\hat Q=\cup_{j\in J(x)}\hat Q_j$
in a neighborhood of $x$ in $x+E$.
 By the inductive assumption
there is an index $\bar j\in J(x)$ and an $x'\in \hat Q_{j_0}$ that belongs to the boundary of $\hat Q_{j_0}$ and $B(x',\ep)\cap\hat Q_j= B(x',\ep)\cap \hat Q$.
for some $\ep>0$.
Finally,
as $0\not\in E$, the same is true for $Q$ and $Q_{j_0}$ as these sets coincide
near $\xb$ with the cones generated by $Q\cap(\xb+E)$ and  $Q_{j_0}\cap(\xb+E)$.\endproof

Returning to the proof of the theorem, set $Q=\im\Phi= \cup_{F\in\cf_C}Q(F)$. If $Q\neq\R^n$, choose a $Q(F)$ and a $\zb\in \bd Q(F)$
according to the lemma, that is such that the intersections of a neighborhood of $\zb$
with $Q(F)$ and $Q$ coincide. In particular the boundaries of the sets near
$\zb$ are identical. The non-singularity condition guarantees that interiors of
all sets $Q(F)$ are nonempty. Therefore, by moving $\zb$ slightly, we can guarantee that $\zb$ belongs to the relative interior of a face of $Q(F)$ of dimension $n-1$.
Clearly, any face of $Q(F)$ is the sum of $F$ and a face of $\Lambda_C(F)$. Thus, there is a $G\in\cf_{\Lambda_C(F)}$ such that $\zb\in\ri(F+G)$ and $\dim(F+G)=n-1$. This may happen
only if $\dim G= n-(\dim F+1)$. By the property (ii) of Definition \ref{dual1}, being a face of $\Lambda_C(F)$, $G$ must be of the form
$G=\Lambda(F')$ for some $F'\in\cf_C$. As $G\subset \Lambda_C(F)$, we necessarily
have $F\subset F'$ and the non-singularity condition implies that $\dim F'=\dim F+1$.

By the assumption, the hyperplane  $E=L(F)+ L(\Lambda_C(F))$ properly
separates $Q(F)$ and $F'$. This however contradicts to the fact that
$\zb$ is a boundary point of $Q(F)$. Indeed, $Q$ and $Q(F)$  coincide in a neighborhood of $\zb$, that is the intersection of $\im \Phi_C$ with the neighborhood lies completely
in one of the  half-spaces determined by $E$. But the relative interior
of $F'$, also a part of $\im\Phi$ lies in the interior of the opposite half-space.
The contradiction completes the proof of the implication and the theorem.
\endproof

The final theorem of the section offers a formula for local regularity rates of $\Phi_C$. As this is going to be a quantitative result, it is more convenient to consider mappings of a somewhat more
general, although formally equivalent, kind, namely
$$
\Phi(x) =  Tx + S(\Lambda_C(F_{min}(x)),
$$
where $T$ and $S$ are linear operators in $\R^n$. We shall do this assuming that $C$ is a cone. This however does not affect the generality of the results as locally
every polyhedral set coincides with a cone (Proposition \ref{loctan}).

\begin{theorem}\label{rate}
Let $K$ be a polyhedral cone and $\Lambda$ a complementarity relation on $K$.
For any $F_1,F_2\in \cf_K$ such that $F_1\subset F_2$ set
$$
r(F_1,F_2)=\inf\{d(T^*z^*,(F_2-F_1)^{\circ}):\; \| z^*\|=1,\; S^*z^*\in (\Lambda(F_1)-\Lambda(F_2))^{\circ}\},
$$
where $T^*$ and $S^*$ stand for the adjoint operators.
Then
$$
\sur \Phi(0|0)= \min\{r(F_1,F_2):\; F_1,\ F_2\in \cf_K, \; F_1\subset F_2  \}.
$$
\end{theorem}

\proof We shall use Theorem \ref{cod} to compute the modulus of surjection of $\Phi$.
To this end we first have to compute tangent cones to the graph of $\Phi$.
Let $z\in \Phi(x)$, that is  $z=Tx+Sy$ for some $y\in F_{max}^H(x)=\Lambda(F_{min}^K(x))$.
Then $(h,w)\in T(\gr\Phi,(x,z))$ is the same as $w=Th+Sv$ and
$(h,v)\in T(\gr F_{max}^H,(x,y))$. So all we need is to compute $T(\gr F_{max}^H,(x,y))$.

Take a pair $(x,y)\in\gr F_{max}^H$ and let  $(h,v)\in T(\gr F_{max}^H,(x,y))$.
We claim that in any case $x+\la h\in F_{max}^K(y)$ for small $\la$.
Indeed, if $H=\{z:\; \lan x_i,z\ran\le 0, \; i=1,\ldots,m \}$, then,
as easily follows from Proposition \ref{minface}
$$
F_{min}^H(y+\la v)=\{z:\; \lan x_i,z\ran=0,\; i\in I(y),\;  \lan x_i,v\ran=0\}.
$$
 But this implies that $F_{min}^H(y)\subset F_{min}^H(y+\la v)$
and therefore $x+\la h\in F_{max}^K(y+\la v)\subset F_{max}^K(y)$ as claimed.

Denote for brevity $F_2= F_{max}^K(y),\ F_1= F_{min}^K(x)$.
As $x$ is in the relative interior of $F_1$, it follows from Proposition \ref{tandif} that
$\cone(F_2-x)= \cone(F_2-F_1)=F_2-F_1$. On the other hand
$\Lambda(F_2)= F_{min}^H(y)$ by definition, so that $y\in{\rm ri}~\!\Lambda(F_2)$
and  $\cone(\Lambda(F_2)-y)= \cone(\Lambda(F_2)-\Lambda(F_2))=L(F_2) $.
Likewise $\cone (F_1-x)=L(F_1)$ and $\cone(\Lambda(F_1)-y)=\cone (\Lambda(F_1)-\Lambda(F_2))$.
Thus, the tangent cone to $\gr F_{max}^K$ at $(x,y)$ contains

\vskip 2mm

\centerline{$\{(h,v):\; h\in L(F_1),\; v\in \Lambda(F_1)-\Lambda(F_2)   \}$}

\noindent and

\centerline{$\{(h,v):\; h\in F_2-F_1,\; v\in L(\Lambda(F_2))   \}$.}

\vskip 1mm

\noindent
The polar of the tangent cone is the same as the polar of its convex hull that
contains the sum of the two cones above.
 It is an easy matter to see (taking consecutively $h=0$ and $w=0$) that
the sum coincides with $(F_2-F_1)\times (\Lambda(F_1)-\Lambda(F_2))$.
We next notice that $F_1\subset F\subset F_2$ for any face $F$ such that $x\in F$ and $y\in \Lambda(F)$.
Therefore $T(F,x)=\cone(F-x)\subset F_2-F_1$ and $T(\Lambda(F),y)=\cone (\Lambda(F)-y)\subset
\Lambda(F_1)-\Lambda(F_2)$. It follows that
 $$
 \conv T(\gr F_{max}^H(x,y))= (F_2-F_1)\times (\Lambda(F_1)-\Lambda(F_2)),
 $$
so that
$$
\conv T(\gr\Phi,(x,z))=\{(h,w):\; h\in F_2-F_1,\; w\in Th + S(\Lambda(F_1)-\Lambda(F_2))\}
$$
Therefore the normal cone to $\gr \Phi$ at $(x,z)$ is the collection of pairs
$(x^*,z^*)$ such that
$$
\lan x^*,h\ran+\lan z^*,Th+Sv\ran\le 0,\quad{\rm if}\; h\in F_2-F_1,\;
v\in \Lambda(F_1)-\Lambda(F_2),
$$
that is
$$
N(\gr\Phi,(x,z))=\{(x^*,z^*): x^*+T^*z^*\in (F_2-F_1)^{\circ},\; S^*z^*\in (\Lambda(F_1)-\Lambda(F_2))^{\circ}   \}
$$
It is an easy matter to see that
$\inf\{\| x^*\|:\; (x^*,z^*)\in N(\gr\Phi,(x,z)),\; \| z^*\|=1  \}$ is precisely
$r(F_1,F_2)$.

To conclude the proof we only need to note that for any pair of faces $F_2,\ F_1\in\cf_K$
such that $F_1\subset F_2$, we can choose $x$ and $y\in F_{max}^H(x)$ such that
$F_2= F_{max}(y),\ F_1=F_{min}^K(x)$, namely any $x\in {\rm ri}~\!F_1$ and
any $y\in {\rm ri}~\!\Lambda(F_2)$ (see Proposition \ref{minface}) and take into account that $\Lambda(F_2)\subset \Lambda(F_1)$). \endproof

\section{Application to variational inequalities.}

In this short section we intend to show that Theorems \ref{main}, \ref{chareg} and \ref{rate}
include Theorems \ref{rob} and \ref{dr96} as particular cases.
Let us start with the proof of the first statement of Theorem \ref{dr96}.
Let $\Phi(x) = Ax+ N(C,x)$.
Let $K$ stand for the tangent cone to $C$ at $\xb$. By Proposition \ref{loctan}
$\Phi_K(x)= Ax + N(K,x)$ is (strongly) regular near $(0,\zb-\xb)$ if
$\Phi$ is (strongly) regular near some $(\xb,\zb)$ and vice versa. Set $\yb=\zb-\xb$.
Then $\Phi_{K^{\circ}}= y+A(N(K^{\circ},y))$ is (strongly) regular near
$(\yb,0)$ if and only if $\Phi_K$ is (strongly) regular near $(0,\yb)$
(Proposition \ref{reg}). Finally, let $K_1=T(K^{\circ},\yb)$. Applying again
Proposition \ref{loctan}, we conclude that $\Phi_1(y)= y+A(N(K_1,y))$
is (strongly) regular near $(0,0)$, hence globally,  if and only if   $\Phi_{K^{\circ}}= y+A(N(K^{\circ},y))$ is (strongly) regular near $(\yb,0)$.

Now, if $\Phi$ is regular near $(\xb,\zb)$, then going down along the above chain of implications, we conclude that $\Phi_1$ is regular near $(0,0)$. But by Theorem \ref{main} this means that $K_1$ is strongly regular near $(0,0)$ and returning up along the same chain of implication we see that $\Phi$ is also strongly regular
near $(\xb,\zb)$. This completes the proof of the first part of Theorem \ref{dr96}.

For the rest  we need to introduce the coherent orientation
and critical face conditions that are central in the theorems stated in Introduction.

\begin{definition}{\rm
It is said that $A$ is
{\it coherently oriented} on $C$ if for each $F\in\cf_C$ the linear operator
defined by
$$
T_F(x) = \left\{\begin{array}{ll} Ax,&{\rm if}\; x\in L(F);\\ x,&{\rm if}\;
x\in L(N(C,F)). \end{array}\right.
$$
is non-singular and  determinants of matrices of all  operators $T_F,\; F\in\cf_C$  have the same sign.

}
\end{definition}

\begin{proposition}\label{coh}
The face separation condition for $(A(C),N(C,\cdot))$ holds if and only if $A$
coherently orientated on $C$.	
\end{proposition}

\proof  We first note that under either of the conditions $A$ must be one-to one on
$L(C)$, so that we may without loss of generality assume that it is one-to-one on the entire space, that is a linear homeomorphism.

The next observation is that to verify the coherent orientation condition, we
only need to compare pairs of operators $T_F$ and $T_{F'}$ such that $F\subset F'$ and $\dim F'=\dim F+1$. Geometrically, this is sufficiently obvious as any two faces
$F,F'\in\cf_C$ can be joined by a chain of faces $F=F_1,\ldots,F_k=F'$ such that
at every link $F_i,F_{i+1}$ one of the faces belongs to the other and
their dimensions differ by one.

On the other hand, the face separation condition for $(A(C),N(C,\cdot))$ means that
for each pair of faces $F,F'$ such that $F\subset F'$ and $\dim F'=\dim F+1$
the subspace $E(F)= L(A(F))+ L(N(C,F'))$ has codimension one and properly separates
$L(A(F))+N(C,F)$ and $A(F')$.

Consider the following  basis $(e_1,\ldots,e_n)$ in $\R^n$:
$(e_1,\ldots,e_{k-1})$ is an orthogonal basis in $L(F)$,
$e_k\in (F'-F)\cap (L(F))^{\perp}$ and  $(e_{k+1},\ldots,e_n)$  is  an orthogonal basis in  $L(N(C,F'))$. 
Set
$$
E_k=\{x=(x_1,\ldots,x_n):\quad x_k=0\},\; E_k^+=\{x:\; x_k>0\}\quad E_k^-=\{x:\; x_k<0\}.
$$
Then $F+ N(C,F')\subset E_k$ and $F+N(C,F)\subset E_k\cup E_k^-$, so that
$$
E=L(A(F))+ L(N(C,F'))=T_F(E_k)
$$
and
$$
\quad A(F)+ N(C,F)=T_F(F+N(C,F))\subset T_F(E_k^-)\cup E.
$$
Thus, the face separation condition means that $Ae_k\in T_F(E_k^+)$.

Let us look now at the coherent orientation condition.
The matrix of $T_F$ in this basis is $(Ae_1,\ldots,Ae_{k-1},e_k,\ldots,e_n)$
and the matrix of $T_{F'}$ is $(Ae_1,\ldots,Ae_k,e_{k+1},\ldots,e_n)$.
Consider the equation $T_{F}x=Ae_k$, and let $x=u+\la e_k$ be the solution, where $u\in E_k$. Clearly,
$$
\la=\frac{{\rm det} T_{F'}}{{\rm det}T_{F}}
$$
The coherent orientation condition means that $\la>0$. This in turn implies that
$x\in E_k^+$ and therefore $Ae_k=T_F(x) \in T_F(E_k^+)$.\endproof

 Theorem \ref{rob} is now immediate. Indeed, by Proposition \ref{comnorm}
 $N(C,\cdot)$ is a complementarity relation on $A(C)$, the proposition along with Theorem \ref{chareg} imply that $\Phi(x)= Ax+ N(C,x)$  is regular if and only if the
 coherent orientation condition is satisfied, and Theorem \ref{main} guarantees that
 regularity and strong regularity are equivalent properties for $\Phi$.

It remains to prove the second part of Theorem \ref{dr96}

 \begin{definition}
 	It is said that {\it critical face condition} is satisfied for $A$ on $C$ if
 	$$
 	z\in\cone(F_2-F_1),\quad A^*z\in (F_2- F_1)^{\circ}\ \Rightarrow\ z=0.
 	$$
 	 for any $F_1,F_2\in\cf_C$,  $F_1\subset F_2$.
 	\end{definition}

 \begin{proposition}\label{lemma}
 	Let $K$ be a polyhedral cone, and let $F, F_2\in\cf_K$, $F_1\subset F_2$. Then
 	$$
 	(F_2-F_1)^{\circ}= N(K,F_1)-N(K,F_2).
 	$$
 \end{proposition}

 \proof By Proposition \ref{tandif} $F_2-F_1= F_2+ L(F_1)$. Therefore
 taking into account that $F_2=K\cap L(F_2)$, we get
 $$
 (F_2-F_1)^{\circ}= (F_2+L(F_1))^{\circ}= F_2^{\circ}\cap (L(F_1))^{\perp}=
 (K^{\circ}+ (L(F_2))^{\perp})\cap L(F_1)= N(K,F_1)+L(F_2)^{\perp}.
 $$
 The last equality needs comments.   First note that the inclusion $\subset$ is obvious. Furthermore,
 $K^{\circ}\cap L(F)^{\perp}=N(K,F)$ for any face $F$ of $K$. Thus 
 $$N(K,F_1)=K^{\circ}\cap L(F_1)^{\perp}\subset (K^{\circ}+ (L(F_2))^{\perp})\cap L(F_1)^{\perp}.$$
 On the other hand, 
 $$(L(F_2))^{\perp}=(L(F_2))^{\circ}\cap (L(F_1))^{\perp}\subset (K^{\circ}+ (L(F_2))^{\perp})\cap L(F_1)^{\perp}$$
 as well.  However $(K^{\circ}+ (L(F_2))^{\perp})\cap L(F_1)$ is a convex cone. Therefore it  contains also  $N(K,F_1)+(L(F_2))^{\perp}$.

  Furthermore,  equalities
 $$
 ({\rm lin}K)^{\circ}=(K\cap(-K))^{\circ}=K^{\circ} +(-K^{\circ})= L(K^{\circ}).
 $$
 hold for any polyhedral cone $K$ and as follows from Proposition \ref{tandif}, $L(F)= {\rm lin}~\!T(K,x)$ for an $x\in{\rm ri }~\!F$. Thus
 $$
 (F_2-F_1)^{\circ}= N(KF_1)+L(N(K,F_2)) = N(K,F_1)-N(K,F_2),
 $$
 the last equality due to the fact that $N(K,F_2)\subset N(K,F_1)$.\endproof

 It follows  that
 $$
 (N(C,F_1)-N(C,F_2))^{\circ} = \cone (F_2-F_1).
 $$
 if $C$ is a polyhedral set and $F_1,F_2\in\cf_C,\; F_1\subset F_2$.
Thus,  the condition
$$
T^*z^*\in (\cone(F_2-F_1))^{\circ},\quad S^*z^*\in (\Lambda(F_1)-\Lambda(F_2))^{\circ} \Rightarrow \ z^*=0
$$
with $T=A,\; S=Id$ and $\Lambda_C=N(C,\cdot)$ reduces precisely to the  critical face condition. Applying Theorem \ref{rate} we get the part (b) of Theorem \ref{dr96}.

\vskip 5mm

\noindent{\bf Some bibliographic comments}.  The interest to piecewise affine mappings
to a large extent associated with the homeomorphism and complementarity problems goes back
as far as to the 50s. The papers by
Dontchev-Rockafellar, Ralph and Robinson quoted in the introduction
contain a fairly detailed account of earlier developments.
Here we just mention a few  results that are most relevant to this study.

First to mention is the paper by Samelson, Thrall and Wesler \cite{STW58}. It seems to be the first to give  a necessary and sufficient condition for the normal map
associated with the standard complementarity problem to be single-valued everywhere.
The condition was stated in terms of positivity of all principal minors of a certain matrix.
Along with that,  a certain separation property is a part of the statement of the theorem in \cite{STW58} and a version of Proposition \ref{coh} for $C=\R_+^n$, the positive orthant in $\R^n$, can be distilled from the result.

Rheinbolt and Vandergraft \cite{RV75} showed that a piecewise mapping is onto if
all associated determinants have the same nonzero sign, and Schramm \cite{RS80},
among other results, proved that under the condition the mapping is open. The inverse implication (openness implies the determinant condition) was established by Eaves and
Rothblum \cite{ER90}.
The last result along with Robinson's principal theorem in \cite{SMR92} (Theorem \ref{rob} here) led Dontchev and Rockafellar in \cite{DR96} to the conclusion that
for any polyhedral $C$
the normal map associated with (\ref{lvi}) as well as
the variational inequality right hand side mapping are locally strongly regular if
they are locally metrically regular.

Finally, we mention the  result of Gowda and Sznajder
\cite{GS95} who found that  any multifunction $\Phi: \R^n\rra \R^n$ whose graph is a union of finitely many polyhedral sets is open if and only if $\Phi^{-1}$ is Lipschitz on the
range of $\Phi$, provided that the range of $\Phi$ is a convex set. As a consequence they
proved that the solution mapping of the linear variational inequality (\ref{lvi}) is nonempty-valued and Lipschitzean on the entire space if and only if
it is single-valued.

\end{document}